\magnification=1200
\overfullrule=0pt
\input amstex
\documentstyle {amsppt}
\parskip=3pt

\def\r{\Bbb R}
\def\c{\Bbb C}
\def\z{\Bbb Z}

\def\CP{\Bbb C \Bbb P}
\def\hh{{\frak h}}
\def\hr{{\hat \rho}}
\def\P{\Bbb P}
\def\PTM{\Bbb P^{TM}}
\def\PE{\Bbb P^E}

\def\so{\frak s\frak o}
\def\hol{\frak h\frak o\frak l}
\def\L{\Lambda}
\def\l{\lambda}
\def\w{\wedge}

\def\<{\langle}
\def\>{\rangle}
\def\nd{\noindent}
\def\CGB{\Cal G_B}
\def\FGB{\frak G_B}
\def\CGE{\Cal G_E}
\def\FGE{\frak G_E}
\def\Hol{\mathop{\hbox{\rm Hol}}\nolimits}
\def\Im{\operatorname{Im}}
\def\im{\operatorname{Im}}
\def\Fr{\operatorname{Fr}}
\def\span{\operatorname{span}}

\topmatter
\title  Parallel Connections Over Symmetric Spaces\endtitle
\author Luis Guijarro, Lorenzo Sadun, and Gerard Walschap \endauthor
\address University of Pennsylvania, Philadelphia, PA 19104\endaddress
\email guijarro\@math.upenn.edu\endemail
\address University of Texas, Austin, TX 78712\endaddress
\email sadun\@math.utexas.edu\endemail
\address University of Oklahoma, Norman, OK 73019\endaddress
\email gwalschap\@math.ou.edu\endemail
\abstract
Let $M$ be a simply connected Riemannian symmetric space, with at most
one flat direction.  We show that every Riemannian (or unitary) vector
bundle with parallel curvature over $M$ is an associated vector bundle
of a canonical principal bundle, with the connection inherited from
the principal bundle.  The problem of finding Riemannian (or unitary) vector
bundles with parallel curvature then reduces to finding
representations of the structure group of the canonical principal
bundle.

\endabstract
\endtopmatter
\document

This paper concerns connections on Riemannian vector bundles over
simply connected symmetric spaces.  Given a hypothesis about the
Riemannian curvature of a symmetric space, we show that every
Riemannian vector bundle with parallel curvature over that space is an
associated vector bundle of a canonical principal bundle, with the
connection inherited from the principal bundle.  The problem of
finding Riemannian vector bundles with parallel curvature then reduces
to finding representations of the structure group of the canonical
principal bundle.  Our results apply also to unitary vector bundles,
since a $\c^k$ bundle can be viewed as an $\r^{2k}$ bundle together
with an additional structure.

The hypothesis we put on a symmetric space is quite mild, and is
satisfied by all simply connected irreducible symmetric spaces, and by
simply connected reducible symmetric spaces that do not contain an
$\r^2$ factor. If a symmetric space does contain such a factor, then
the theorem does not hold, and we give explicit counterexamples.

Recall that a Riemannian connection on a vector bundle $E\to M$ is said to be
{\sl Yang-Mills\/} if its curvature tensor $R^E$ is harmonic; i.e., if
$d_{A}^* R^E=0$, where $d_{A}^*$ is the covariant divergence
operator. Yang-Mills
connections have generated substantial interest, and much of our current
knowledge of
the topology of smooth 4-manifolds comes from their study [DK]. In the special
case of the Levi-Civita connection on the tangent bundle of a Riemannian
manifold, there is a stronger concept which has been of central importance in
Riemannian geometry. A Riemannian manifold $M$, with Riemann curvature tensor
$R^M$, is said to be {\sl locally symmetric\/} if all covariant derivatives
of $R^M$ vanish. This is equivalent to $R^M(X,Y)Z$ being a parallel vector
field along every path $\gamma$
for any $X$, $Y$, $Z$ parallel along $\gamma$.
This concept generalizes naturally to arbitrary
Riemannian vector bundles $E\to M$:

\definition{Definition}
A connection on $E$ is said to be {\sl
parallel} if $\nabla R^E=0$; equivalently,
if for any smooth path $\gamma$ of $M$, parallel vector fields $X$, $Y$
along $\gamma$, and parallel section $U$ of $E$ along $\gamma$,
the section $R^E(X,Y)U$ is parallel along $\gamma$.
\enddefinition

This condition involves, of course, the Riemannian metric on the base
manifold $M$.  (By a Riemannian vector bundle we mean a vector bundle
with an inner product on each fiber.  A connection on such a bundle is
required to respect the inner product, meaning that the inner product
is parallel. In this paper all vector bundles are assumed to be real
and Riemannian, except where stated otherwise.)

Recall that a symmetric space is a Riemannian manifold $M$ for which
the geodesic symmetry at any point is a global isometry of $M$. Every
locally symmetric space has a globally symmetric universal
cover. Conversely, every symmetric space is locally symmetric.  If $M$
is a symmetric space, then $M$ can be written as $G_0/H_0$, where
$G_0$ and $H_0$ are groups, with $H_0=\Hol(M)$, the holonomy group of
$M$.  For example, $S^n=SO(n+1)/SO(n)$ and $\CP^n =
[U(n+1)/U(1)]/U(n)$.  By lifting to the universal cover of $G_0$ we can
write $M=G/H$, where $G$ is simply connected and $H$ is a cover of
$H_0$.  If $M$ is simply connected, then $H$ is connected. For
example, $S^n = Spin(n+1)/Spin(n)$ and $\CP^n = SU(n+1)/U(n)$. (In the
$\CP^n$ example, $H_0$ and $H$ are both isomorphic to $U(n)$ as
groups, but $H$ is still the $n+1$-fold cover of $H_0$.)

Let $M_p$ denote the tangent space of $M$ at $p$. The space $\Lambda^2
(M_p)$ is naturally identified with $\so(M_p)$, the Lie algebra of
skew-adjoint endo\-morphisms of $M_p$, by requiring that
$$u\w v\,(w)=\<u,w\>v-\<v,w\>u, \qquad u,v,w\in M_p . \tag{1.1}$$
The holonomy of $M$ acts on $\L^2$ in a natural way via $h(u\w v)=hu\w
hv$, for $h\in\Hol(M)$. Under the above identification, the action of
$\Hol$ is by
conjugation on $\so(M_p)$. This is also true for $\Lambda^2E$ and the
holonomy group in any vector bundle $E$ with a Riemannian
connection. Furthermore, the 
Riemannian curvature $R^M$ at $p$ is an equivariant map from
$\so(M_p)$ to itself. It is easy to see that in general $[\ker R^M,
\Im R^M]\subset\ker R^M$, and therefore the same is valid for the
linear subspace spanned by $[\ker R^M,
\Im R^M]$.  For our purposes, we will consistently
assume that $R^M$ satisfies a stronger condition:

\nd {\bf Condition A}: $\span[\ker R^M, \Im R^M] = \ker R^M$.

This condition is well understood; in fact,
in \S 4 we will prove:
\proclaim{Theorem 1.1} Let $M$ be a simply connected Riemannian symmetric
space.  Condition A holds unless $M$ is the product of\/ $\r^2$ and
another symmetric space, in which case Condition A fails.
\endproclaim

Put another way, Condition A holds as long as the dimension of the
center of $G$ does not exceed 1. This is discussed in detail in
section 4.

Since the Levi-Civita connection on the tangent bundle of a symmetric
space $M$ is parallel, any bundle built from the tangent bundle will
also admit a parallel connection, naturally induced from the
Levi-Civita connection on $TM$.  In particular ${\otimes \atop k} TM$
and $\Lambda^k(TM)$ admit parallel connections, as does the
(principal) bundle $\Fr(TM)$ of orthonormal frames of $TM$.
Sub-bundles of $\Fr(TM)$, and finite covers of $\Fr(TM)$ or of
sub-bundles of $Fr(TM)$, also inherit parallel
connections, whenever such subbundles are invariant by parallel transport.

If $B$ is a principal $H$-bundle over $M$ and $\rho: H \to SO(k)$
is a representation of $H$, then the associated vector bundle $E :=B
\times_\rho \r^k$ is the quotient of $B \times \r^k$ by the
equivalence $(bh,v) \sim (b,\rho(h)v)$. The equivalence class of
$(b,v)$ is denoted $[(b,v)]$.  The vector bundle $E$ inherits a
connection from that of $B$, which will be parallel if the connection
on $B$ is parallel.  Thus vector bundles with parallel connections
naturally arise from representations of the holonomy group of a
symmetric space $M$, or of its finite covers.

Parallel connections are also related to group actions.  Let $G$ be a
group that acts orthogonally and transitively on a vector bundle $E$,
covering a group  of isometries of $M$, and suppose that $E$ has a
connection that is 
invariant under this $G$-action.  Then that connection is known to be
parallel.  Furthermore, for any fixed bundle $E$ and $G$-action,
Wang's theorem [W] classifies the $G$-invariant connections.

In this paper we show that, on symmetric spaces satisfying Condition A,
these two constructions are equivalent, and are the only source of parallel
connections:

\proclaim{Main Theorem} Let $M=G/H$ be a simply connected symmetric space
written as a canonical group quotient, with
$G$ simply connected, and suppose that the canonical
metric on $M$ satisfies 
Condition A.  Let $E$ be any rank-$k$ real vector bundle with parallel
connection over $M$.   Then there exists a representation $\rho: H \to SO(k)$
such that $E$ is isomorphic to the vector bundle $E'=G \times_\rho \r^k$, with
the isomorphism taking the natural connection on $E'$ to the given
connection on
$E$.
\endproclaim

\proclaim{Corollary 1} Every parallel connection over $M$  is invariant under a
$G$-action on the bundle that covers the natural $G$-action on $M$.
\endproclaim

\proclaim{Corollary 2} Every unitary rank-$k$ vector bundle with parallel
connection over $M$ is isomorphic to $G \times_{\rho'} \c^k$, where
$\rho': H \to U(k)$ is a rank-$k$ unitary representation of $H$.
\endproclaim

Several authors have investigated the question of when the action of a
given transitive group of isometries on the base $M$ can be lifted to
a bundle over $M$, see for instance \cite{B-H}. Since such a lifted
action implies the existence of a parallel connection, Corollary 1
provides a complete answer: Over a symmetric space satisfying the
hypotheses of the main theorem, the vector bundles that admit lifts of
transitive group actions are exactly the bundles that admit parallel
connections, which are exactly the associated vector bundles of
$G$ itself. Notice, though, that $G$ is usually a cover of the isometry group
$G_0$ of the base, and that the action of the latter will not, in general, lift
to the bundle: For example, the action of $SO(5)$ on the 4-sphere does not lift
to the spinor bundle, even though that of $Spin(5)$ does.

The main theorem tells us that a local condition on a bundle with
connection, namely that its curvature be parallel, implies a rigid
global structure.  However, for this to be true we need assumptions
about the underlying manifold, namely that it be a simply connected
symmetric space whose metric connection satisfies Condition A. Without
these conditions the conclusions of the theorem can be false, as the
following counterexamples show.

Let $M=T^2=G/H$, where $G=\r^2$ and $H=\z^2$. $M$ is not simply
connected and Condition A fails, as $\ker R^M$ is all of $\Lambda^2
(M_p)$ while $\Im R^M$ is zero.  Since $R^M=0$, any bundle constructed
from a representation of $H$ is necessarily flat (although possibly
with nontrivial holonomy).  However, there exist bundles and
connections over $T^2$ with parallel nonzero curvature.  For example,
consider the trivial complex line bundle over $\r^2$.  Taking the
quotient of this by the relation $f(x+n,y+m)= \exp(2\pi i n y) f(x,y)$
gives a complex line bundle of Chern class $+1$ over $T^2$ that we
denote $E$. $E$ admits many connections of constant nonzero curvature;
one has connection form $A= -2 \pi i x dy$. Not only does this line
bundle not come from a representation of $H$, but the connection is
not invariant under translation.  Indeed, the curvature itself
measures the extent to which holonomy changes when a homologically
nontrivial loop is translated.  There is a 2-parameter family of
parallel connections on $E$, indexed by $T^2$ itself, and translation
takes one such connection into another.

Now let $M=G=\r^2$, with $H$ trivial.  $M$ is simply connected but
does not satisfy Condition A.  Once again, any bundle built
from the tangent bundle with its canonical structure has a flat
connection, but there are connections over $\r^2$ that are parallel
but not flat.  For example, one can again look at a complex line
bundle with connection form $A=-2 \pi i x dy$ and constant curvature
$-2 \pi i dx \wedge dy$.  Unlike the $T^2$ example, this connection is
invariant under translation, so Wang's theorem does apply.

The proof of the main theorem is based on the following
observations. Let $M=G/H$ be a simply connected symmetric space, with
$G$ simply connected.  $H$ is a covering space of $\Hol(M)$, the
holonomy group of $M$.  Let $\pi:H \to \Hol(M)$ be the covering map,
so that $TM=G \times_\pi \r^n$.  Suppose $E= G \times_\rho R^k$, for a
representation $\rho: H \to SO(k)$, and endow $E$ with the connection
induced by the one on $G$. Then the holonomy around a loop $\gamma$ 
for the three bundles $TM$, $G$, and $E$ are related by the
commutative diagram
$$\matrix
\gamma\cr
\noalign{\vskip6pt}
\llap{$\scriptstyle \P_\gamma^M$} \swarrow\quad
\big\downarrow \rlap{$\scriptstyle\P_\gamma^G$}
\quad \searrow\rlap{$\scriptstyle\P_\gamma^E$}\cr
\noalign{\vskip6pt}
\Hol (M_p)
\buildrel \pi\over\longleftarrow
\enspace H \enspace
\buildrel \rho\over\longrightarrow
\Hol (E_p),
\endmatrix\tag{1.2}$$
where $\PE_\gamma$ denotes parallel translation in $E$ along the curve
$\gamma$ in $M$, etc.  One cannot compute $\PE_\gamma$ directly from
$\PTM_\gamma$, since $\pi$ is not 1--1.  However, looking at
infinitesimal loops we have the commutative diagram
$$\matrix
\Lambda^2 (M_p)\cr
\noalign{\vskip6pt}
\llap{$\scriptstyle R^M$} \swarrow\quad
\big\downarrow \rlap{$\scriptstyle R^G$}
\quad \searrow\rlap{$\scriptstyle R^E$}\cr
\noalign{\vskip6pt}
\hol (M_p)
\buildrel {\hat\pi}\over\longleftarrow
\enspace \hh \enspace
\buildrel {\hat\rho}\over\longrightarrow
\hol (E_p),
\endmatrix\tag{1.3}$$
where $\hat \pi: \hh \to \so(M_p)$ and $\hat \rho: \hh \to \so(E_p)$
are the Lie algebra homomorphisms corresponding to the group
homomorphisms $\pi$ and $\rho$.  Since $\pi$ is a covering map, $\hat
\pi$ is an algebra isomorphism, and we have
$$
R^E = \hat \rho \circ (\hat \pi)^{-1} \circ R^M. \tag{1.4}
$$
Since $R^M$ is parallel and $\hat \rho$ and $\hat \pi$ are fixed,
$R^E$ is parallel.  Notice that $\ker R^M \subset \ker R^E$, so that
$R^E \circ (R^M)^{-1}$ is well-defined (although $(R^M)^{-1}$ may not
be), and is a Lie algebra homomorphism from $\hol(M_p)$ to $\hol(E_p)$.
Indeed, it is equal to $\hat \rho \circ \hat \pi^{-1}$.

Now suppose that $E$ is a vector bundle over $M=G/H$ with parallel
curvature.  To show that $E$ is an associated vector bundle of $G$, we
reverse the above argument.  First we show that $\ker R^M \subset \ker
R^E$. This is the step where condition A is needed. Then, with
$\hh$ identified with $\hol(M_p)$, we define $\hat 
\rho = R^E \circ (R^M)^{-1}$ and show that it is a Lie algebra
homomorphism from $\hh=\hol(M_p)$ to $\hol(E_p)$.  Next, we
exponentiate this algebra homomorphism to get a group homomorphism
$\rho: H \to \Hol(E_p) \subset SO(k)$.  This is somewhat delicate,
since $H$ may not be simply connected.  To get around this, we extend
$\hat \rho$ to an algebra homomorphism from $\frak g$ to the Lie
algebra of the extended gauge group, exponentiate that to a group
homomorphism from $G$ to the extended gauge group, and then restrict
the group homomorphism to $H$. This is the step that requires the
simple connectivity of $G$.

Having reconstructed $\rho$ from the curvature of $E$, we define
$E' = G \times_\rho \r^k$.  Identifying the fibers of $E$ and $E'$ at our base
point $p$, we note that the curvatures of $E$ and $E'$ are equal, and
indeed the
two vector bundles have the same holonomy around any given loop.  We then use
parallel transport to identify the fibers of
$E$ and $E'$ at every point, completing the proof.

In \S 2 we establish the local properties of $R^M$ and $R^E$, where
$E$ is any vector bundle with parallel curvature.  In particular, we show
(Theorem 2.4) that $R^E \circ (R^M)^{-1}$ is a Lie algebra homomorphism. In
\S 3
we construct the bundle $E'$ and show that it is isomorphic to
$E$.  In \S 4 we consider Condition A, and show that it holds on all  simply
connected symmetric spaces that do not have an $\r^2$ factor. Finally, in \S 5
we discuss some of the implications of the Main Theorem.

\head 2.\quad Local properties of parallel connections \endhead

We briefly recall some basic facts about holonomy and curvature in
bundles. The reader is referred to \cite{P} for further details and
other facts that will be freely used here. It will often be convenient
to work with principal bundles instead of vector bundles. Recall that
an oriented Riemannian vector bundle $E$ of rank $k$ over $M$ is
isomorphic to $\Fr(E)\times_{\rho_0}\r^k$, where $\Fr(E)$ is the
principal $SO(k)$ bundle of oriented orthonormal frames of $E$ and
$\rho_0$ is the fundamental representation of $SO(k)$.  A sub-bundle
of $\Fr(E)$ is obtained by taking all parallel translates of a given
frame.  This is a principal bundle over $M$, whose fiber is the
holonomy group $\Hol(E)$ of $E$.

A connection on $E$ induces one on the corresponding principal bundle
$\Fr(E)$, on all sub-bundles of $\Fr(E)$, and on all covers of
sub-bundles of $\Fr(E)$. Conversely, if $B$ is a principal $H$-bundle
over $M$ and $\rho$ is an orthogonal $k$ dimensional representation of
$H$, then a connection on $B$ induces one on the vector bundle
$E=B\times_\rho\r^k$, which will be parallel if the connection on $B$
is parallel.

The curvature tensor $R^E$ of a connection on $E$ is a bundle map 
$\L^2TM\to\hol(E)\subset \so(E)$. A result of Ambrose and Singer says
that the fiber of $\hol(E)$ over $p\in M$ is spanned by the parallel translates
of all curvature transformations $R^E(x\w y)$, $x\w y\in\L^2M_q$, $q\in M$,
along smooth curves from $q$ to $p$ \cite{P}. When the connection is
parallel, the 
curvature is invariant under parallel translation, and we obtain:

\proclaim{Lemma 2.1} Let $E$ be a vector bundle with parallel connection
over a manifold $M$, and let $p$ be a point on $M$, then
\roster \item $R^E:\L^2(M_p)\to\hol(E_p)$ is onto.
\item If $\gamma$ is a curve in $M$ with $\gamma(0)=\gamma(1) =p$, and $A\in
\L^2(M_p)=\so(M_p)$, then
$$R^E\circ\P_\gamma^{\so(TM)}(A)=R^E(\PTM_\gamma A\, \PTM_{-\gamma}) =
\PE_\gamma R^E(A)\,\PE_{-\gamma}=\P_\gamma^{\so(E)}\circ R^E(A).$$
\item For $A,B\in\L^2(M_p)$, $R^E[R^MA,B]= [R^EA,R^EB]$.
\endroster
\endproclaim

\demo{Proof} \therosteritem1 and \therosteritem2 are clear from the above
discussion, together with the fact that the holonomy acts by
conjugation on both $\so(TM)$ and $\so(E)$. It suffices to establish
\therosteritem3 for decomposable elements
$A=x\w y$, $B=z\w w$, $x,y,z,w\in M_p$. So consider vector fields $X$,
$Y$ on $M$ such that $X_p=x$, $Y_p=y$, and $[X,Y]=0$ in a neighborhood
of $p$. If $\{\phi_s\}$ and $\{\psi_s\}$ denote the local 1-parameter
groups of $X$ and $Y$, then for small $t$, the product $\gamma_t$ of
the curves $s\mapsto\phi_{\sqrt s}(p)$, $s\mapsto\psi_{\sqrt s}
\phi_{\sqrt t}(p)$, $s\mapsto\phi_{-\sqrt s}\psi_{\sqrt t}\phi_{\sqrt t}(p)$,
and $s\mapsto\psi_{-\sqrt s}\phi_{-\sqrt t}\psi_{\sqrt t}\phi_{\sqrt t}(p)$,
$0\le s\le t$, is a piecewise smooth loop at $p$, such that
$$
\frac d{dt}\Big |_0\PTM_{\gamma_t}=R^M(A),\qquad \text{ and }
\frac d{dt}\Big |_0\PE_{\gamma_t}=R^E(A),\tag{2.2}
$$
cf\. also \cite{P}. If $\beta$ denotes the curve $t\mapsto
R^E\P^{\so(TM)}_{\gamma_t}(B)$ in $\L^2(M_p)$, we have that
$\beta(t)=R^E\PTM_{\gamma_t}B(\PTM_{\gamma_t})^{-1}$, and by the first identity
in \thetag{2.2}, $\dot\beta(0)=R^E[R^MA,B]$. On the other hand, $\beta(t) =
\P^{\so(E)}_{\gamma_t}R^E(B)=\PE_{\gamma_t}R^E(B)(\PE_{\gamma_t})^{-1}$,
and the
second identity in \thetag{2.2} yields $\dot\beta(0)=[R^E(A),R^E(B)]$, which
establishes the claim. \qed
\enddemo

We now consider the case when the base is a locally symmetric space. The
curvature tensor of the Levi-Civita connection then satisfies the following
property:

\proclaim{Proposition 2.2} Let $M$ be a locally symmetric space,
$E_\l\subset\so(M_p)$ the $\l$-eigenspace of the curvature operator
$R^M:\so(M_p)\to\so(M_p)$. If $\l\ne0$, then $E_\l$ is an ideal of the
algebra bundle $\hol(M_p)\subset\so(M_p)$. Thus, $\hol(M_p)$ is a Whitney
sum $\oplus E_i$ of parallel eigenspace bundles, which are fiberwise
ideals.
\endproclaim

\remark{Remark} Adding the 0-eigenbundle of the curvature decomposes $\so(TM)$
as a Whitney sum of parallel subbundles. In general, however, the kernel of
$R^M$ is not closed under the Lie bracket.
\endremark

\demo{Proof of Proposition 2.2}
Let $A$ and $B$ denote $\l$-eigenvectors. Then Lemma
2.1\therosteritem3 with $E=TM$ implies
$$\l R^M[A,B]=R^M[R^MA,B]=[R^MA,R^MB]=\l^2[A,B],
\tag{2.3}
$$ so that $E_\l$ is a subbalgebra bundle if $\l\ne0$. Similarly, for nonzero
$\mu\ne\l$, and $\l$-eigenvector $A$, $\mu$-eigenvector $B$, we have
$$
\eqalign {\l R^M[A,B]=R^M[R^MA,B]= & [R^MA,R^MB] \cr = &
R^M[A,R^MB]=\mu R^M[A,B],}
\tag{2.4}
$$ where the identity $[R^MA,R^MB]=R^M[A,R^MB]$ follows by
interchanging $A$ and $B$ in Lemma 2.1\therosteritem3. Thus, $[A,B]\in\ker
R^M\cap \hol(M_p)=\{0\}$. \qed
\enddemo

\proclaim{Lemma 2.3} If $M$ is a locally symmetric space satisfying
Condition A and $R^E$ is the curvature of a parallel connection on a
bundle $E$ over $M$, then $\ker R^M \subset \ker R^E$.
\endproclaim

\demo{Proof} Suppose that $A \in E_\l$ and that $B \in \ker R^M$, with
$\lambda \ne 0$.  Then
$$
\l R^E[A,B]=R^E[R^MA,B]=  [R^EA,R^EB] =  R^E[A,R^MB]=0.
\tag{2.5}
$$ Thus $R^E[E_\l,\ker R^M]=0$.  Summing over $\l$ we have that
$R^E[\Im R^M, \ker R^M]=0$.  But by Condition A,
$\span[\Im R^M,\ker R^M]=\ker R^M$, so $R^E(\ker R^M)=0$. \qed
\enddemo

Since $R^M$ is onto $\hol(M_p)$, and  $R^E(\ker R^M)=0$,  there exists
a map $\hat \rho: \hol(M_p) \to \so(E_p)$ such that
$R^E = \hat \rho \circ R^M$.  Formally we write $\hat \rho = R^E
\circ (R^M)^{-1}$.

\proclaim{Theorem 2.4} Let $M$ denote a symmetric space,  $R^E$ the
curvature at a point $p$ of
a parallel connection on a bundle $E$ over $M$. Then $\hat\rho:=R^E\circ
(R^M)^{-1}:\hol(M_p)\to\so(E_p)$ is a Lie algebra bundle homomorphism whose
image is all of $\hol(E_p)$.
\endproclaim

\demo{Proof} Since $R^E = \hat \rho \circ R^M$, the image of $\hr$ is the
image of $R^E$, namely $\hol(E_p)$.  To see that $\hr$ is a Lie
algebra homomorphism, we need to check that $\hr[A,B]=[\hr(A),\hr(B)]$
when $A$ and $B$ belong to a common eigenspace $E_\l$ of $R^M$, and
when $A\in E_\mu$, $B\in E_\l$, $\l\ne\mu$. Let $b = \tfrac1\l B\in
(R^M)^{-1}B$.  In the former case,
$$
[\hr(A),\hr(B)]=R^E[A,b]=\tfrac1\l R^E[A,B]=\hr[A,B],
\tag{2.6}
$$
since $[A,B]$ belongs to $E_\l$. In the latter case, we have
$$
[\hr(A),\hr(B)]=R^E[A,b]=\tfrac1\l R^E[A,B].
\tag{2.7}
$$
After interchanging the roles of $A$ and $B$, we see that the left
side of (2.7) is also equal to $(1/\mu) R^E[A,B]$, and therefore
vanishes.  However, $[A,B]=0$ by Proposition 2.2, so
$[\hr(A),\hr(B)]=0 =\hr[A,B]$.  \qed
\enddemo

\head 3.\quad The global structure of $E$ \endhead

In this section we complete the proof of the main theorem, using the
local structure of \S 2 to construct a bundle $E'$, and then showing
that $E$ and $E'$ are isomorphic with the isomorphism respecting the
connections.   Our base manifold is a Riemannian symmetric
space $M$. Recall that $M=G_0/H_0$, where $G_0$ is the largest connected group
of isometries of $M$, and $H_0$ is the isotropy group at some fixed point $p$.
If $\lambda:H_0\to SO(M_p)$ denotes the (faithful) linear isotropy
representation $\lambda(h)=h_{\ast p}$, then the tangent bundle $TM$ is
isomorphic to $G_0\times_\lambda M_p$ via $[(g,u)]\mapsto g_{\ast p}u$, and
the holonomy group of $M$ at $p$ is $\lambda(H_0)$. An equivalent description
of the tangent bundle is obtained by considering the involutive automorphism
$\sigma$ of $G_0$ given by $\sigma(g)= s_p\circ g\circ s_p^{-1}$, where $s_p$
denotes the geodesic symmetry at $p$. Its derivative at the identity squares
to 1 and decomposes the Lie algebra $\frak g$ into a sum $\frak g=\frak
h\oplus\frak m$ of the plus and minus 1 eigenspaces of $\sigma_{*e}$. The
derivative of the projection $\pi:G_0\to M$ (which sends $g$ to $g(p)$) has
kernel $\frak h$ and thus identifies $\frak m$ with the tangent space of $M$
at $p$. A straightforward computation shows that for $h\in H_0$, the diagram
$$\CD
\frak g		@>Ad_h>> \frak g \\
@V\pi_{*e}VV			      @VV\pi_{*e}V \\
M_p	    @>h_{*p}>>    M_p
\endCD \tag{3.1}
$$
commutes, so that $TM$ is also identified with $G_0\times_{Ad}\frak m$, with
$Ad$ denoting the adjoint representation of $H_0$.

 We
work throughout with
$M=G/H$ where $G$ is the universal cover of $G_0$,
and identify the Lie algebra $\hh$ of $H$ with $\hol(M_p)$.  Our first
task is to promote the Lie algebra
homomorphism $\hr: \hol(M_p) \to \so(E_p)$ into a group homomorphism
$\rho: H \to SO(E_p)$.  Since in general $H$ is not simply connected,
this promotion is not automatic.  We must extend $\hr$ to a
homomorphism from $\frak g$ to an appropriate Lie algebra, promote
that extension to a group homomorphism, and then restrict the
homomorphism to $H$.

Recall that if $B$ is a principal $K$-bundle over $M$, the {\sl gauge
group\/} $\CGB$ of $B$ is the group of $K$-equivariant diffeomorphisms
of $B$ that cover the identity on $M$, cf\. e.g\. \cite{BL}. It is
canonically identified with the space of sections of the bundle of
groups $B\times_c K$, where $c$ denotes conjugation; i.e.,
$(bk_1,k_2)\sim(b, k_1k_2k_1^{-1})$. The gauge algebra $\FGB$ of $B$
is similarly defined as the Lie algebra of sections of the Lie algebra
bundle $B\times_{\operatorname{Ad}}\frak K$, and there is a fiberwise
exponential map
$$
\exp:\FGB\to\CGB. \tag{3.2}
$$
In the case of a Riemannian vector bundle $E$, the gauge group $\CGE$
and algebra $\FGE$ are defined to be the gauge group and algebra of
the principal bundle $\Fr(E)$ of oriented orthonormal
frames. Equivalently, $\CGE$ is the group of orthogonal bundle maps $E \to E$
that cover the identity on $M$.  An element of $\frak G_E$ can also
be viewed as a (vertical) vector field on $E$.

Let $M=G/H$ be a symmetric space with $G$ simply connected, and $E$ a
vector bundle over $M$. Define the {\sl enlarged gauge group\/} of $E$
to be the group $\tilde\Cal G_E$ of all orthogonal bundle maps $E\to
E$ that cover the natural action of $G$ on $M$. The projection
$\alpha: E \to M$ determines an exact sequence
$$
1\to\Cal G_E\to\tilde\Cal G_E\overset\alpha\to\longrightarrow G\to1,
\tag{3.3}
$$
which, when differentiated at the identity, yields the exact sequence
of vector spaces
$$
0\to\frak G_E\to\tilde\frak G_E\overset{\alpha_*}\to\longrightarrow\frak
g\to0. \tag{3.4}
$$
A connection on $E$ determines a splitting of the sequence (3.4). An
element $X\in\tilde\frak G_E$, viewed as a vector field on $E$, splits
accordingly as $X=X^h+X^v$, where $X^v$ is a section of $\hol(E)$,
and $\alpha_*X^h=\alpha_*X\in\frak g$.

When $E=TM$, any $X\in\frak g$ induces a vector field $\tilde
X\in\tilde\frak G_M$ on $TM$: Identify $X$ with the Killing field on
$M$ induced by the $G$-action. If $\phi_t$ denotes the one-parameter
group generated by $X$, then $\phi_t{}_*$ is a bundle map that covers
$\phi_t$, and we set $\tilde X:=\left.{\frac d{dt}}\right|_0(\phi_t{}_*)$. By
standard properties of connections (see e.g\. \cite{P}), one has
$$
[\tilde X^v,\tilde Y^h]=0,\quad [\tilde X^h,\tilde Y^h]^v=R^M(X,Y),\qquad
X,Y\in\frak g. \tag{3.5}
$$
Define a map $\phi:\frak g\to\tilde\frak G_E$ by $\phi(X)=  X^h+
\hr(\tilde X^v)$, where $\hr=R^E\circ( R^M)^{-1}: \hol(TM)\to\hol(E)$
is the Lie algebra bundle homomorphism from Theorem 2.4, and $ X^h$ is
the horizontal lift of $ X$ (viewed as a Killing field on $M$) to $E$.

\proclaim{Lemma 3.1} $\phi:\frak g\to\tilde\frak G_E$ is a Lie algebra
homomorphism.
\endproclaim
\demo{Proof} We compute
$$\split \widetilde{[X,Y]} &= [\tilde X,\tilde Y]= [\tilde X, \tilde Y]^h +
[\tilde X, \tilde Y]^v= [\tilde X^h, \tilde Y^h]^h + [\tilde X,\tilde Y]^v\\ &=
[\tilde X^h, \tilde Y^h]^h + [\tilde X^h, \tilde Y^h]^v+[\tilde X^v,\tilde
Y^v]\\&=[\tilde X^h, \tilde Y^h]^h + R^M( X\wedge  Y) + [\tilde
X^v,\tilde Y^v].
\endsplit \tag{3.6}
$$ Thus,
$$ \split \phi[X,Y] &  = [ X^h, Y^h]^h+ R^E( X\wedge  Y) +
\hr[\tilde X^v, \tilde Y^v]\\ & = [ X^h, Y^h]^h+ R^E( X\wedge  Y) +
[\hr(\tilde X^v), \hr(\tilde Y^v)],
\endsplit \tag{3.7}
$$ since $\hr$ is an algebra homomorphism. On the other hand,
$$\split [\phi X,\phi Y]&=[  X^h+\hr(\tilde X^v),   Y^h+\hr(\tilde Y^v)]\\&=
[ X^h, Y^h]^h +[ X^h, Y^h]^v + [\hr(\tilde X^v), \hr(\tilde Y^v)]\\ &=[ X^h,
Y^h]^h+R^E(X\wedge Y) + [\hr(\tilde X^v), \hr(\tilde Y^v)]. \qed
\endsplit \tag{3.8}
$$
\enddemo

\proclaim{Proposition 3.2}
Let $p$ denote the coset $eH\in G/H=M$. There exists a group
homomorphism $\rho: H \to SO(E_p)$ whose derivative at the identity is
$\hr|_p: \hh \to \so(E_p)$.
\endproclaim

\demo{Proof}  Since $G$ is simply connected, the algebra homomorphism
$\phi$ can be promoted to a group homomorphism $\Phi: G \to \tilde
\CGE$.  Restricting $\Phi$ to $H$ we get a map $H \to \tilde \CGE$.
Now, if $X \in \frak h$, then $ X$, viewed as a Killing field on $M$,
has a zero at $p$ because $H_0$ is the isotropy group of $p$.
Therefore $\Phi(\exp(tX))$, and generally $\Phi(h)$ for any $h \in H$,
preserves the fiber $E_p$.  We define $\rho(h)$ to be the restriction
of $\Phi(h)$ to that fiber.  $\rho$ is a group homomorphism since
$\Phi$ is, and the derivative of $\rho$ is $ev_p \circ \phi$, where
$ev_p$ is evaluation at the base point $p$.  But for $X \in \frak h$,
$$ev_p \circ \phi(X)=\tilde X^h(p) + \hr(\tilde X^v)(p)=0+\hr_p(X).
\qed \tag{3.9}
$$
\enddemo

We now define the bundle $E'$ to be $G \times_\rho E_p$.  $E'$ inherits a
natural parallel connection from $G$. The Main Theorem boils down to

\proclaim{Theorem 3.3} $E$ and $E'$ are isomorphic vector bundles, with the
isomorphism preserving the connection.
\endproclaim

\demo{Proof} The fibers of
$E$ and $E'$ over $p$ are canonically identified.  Denote the
isomorphism $I_p:E'_p \to E_p$.  $I_p$ is an element of
$Hom(E',E)$. We must show that $I_p$ can be extended in a parallel
fashion to all of $M$.

Let $F$ be the vector bundle $Hom(E',E)$. $F$ inherits a parallel
structure from those of $E'$ and $E$.  In particular, the holonomy of
$F$ is generated by the curvature of $F$ at $p$, and this curvature is
constructed from the curvatures of $E'$ and $E$.  By (1.4) the
curvature of $E'$ at $p$ is
$$ R^{E'}_p = \hr_p \circ R^M_p = R^E_p. \tag{3.10}
$$
Thus, for any $A \in \Lambda^2(M_p)$, $R^{E'}(A)=R^E(A)$, and so
$$ R^F(A)(I_p) = 0. \tag{3.11}
$$
Since $I_p$ is in the kernel of $R^F(A)$ for every $A$, $\Hol(F)$ acts
trivially on $I_p$, so $I_p$ can be extended, by parallel transport,
to a parallel section of $F$.  This section is our global isomorphism
between $E'$ and $E$. \qed
\enddemo

This completes the proof of the main theorem.  Corollary 1 is an
obvious consequence.  For Corollary 2, note that a rank-$k$ complex
vector bundle is just a real rank-$2k$ vector bundle with an
additional feature, namely the complex structure $J$.  If $E$ is a
unitary bundle with parallel connection, then the holonomy of $E$
commutes with $J$, and so $\hr_p$ commutes with $J$.  Thus $\rho$ is
actually a rank-$k$ unitary representation of $H$, not just a rank
$2k$ orthogonal representation.

In the proof of the Main Theorem we have constructed the
representation $\rho$ as linear operators on the fixed vector space
$E_p$ rather than as matrices acting on $\r^k$, and the definition of
$E'$ is similar.  The difference between working with $E_p$ and $\r^k$
is of course nothing more than a choice of orthonormal basis for
$E_p$.

\head 4.\quad The curvature condition \endhead

In this section we prove Theorem 1.1, showing that Condition A holds  on all
simply connected symmetric spaces that do not have an $\r^2$ factor. We begin
with a special case:

\proclaim{Theorem 4.1} Condition A applies on all simply connected irreducible
Riemannian symmetric spaces.
\endproclaim

\demo{Proof}  Let $M$ be such a space of dimension $n$, let $p\in M$ be a
point, and suppose that $W := \span[\ker R^M_p, \Im R^M_p]\subset
\so(M_p)$ is a proper subspace of $V := \ker R^M_p$. Let $A$ be an
element of $V$, orthogonal to $W$.  $A$ is then orthogonal to all of
$[\so(M_p), \Im R^M]$, since $[\Im R^M, \Im R^M]$ is contained in
$\Im R^M$, and so is orthogonal to $A$.  Since $\Hol$ acts
orthogonally on $\so(M_p)$, $A$ is invariant under $\Hol$.

Since $A$ is a skew endomorphism of $M_p$, it can be expressed as a
block diagonal $n \times n$ matrix, where each block is either zero or
a skew $2 \!\times\! 2$ matrix $L_\lambda=\left ( \matrix 0 & \lambda
\\ -\lambda & 0 \endmatrix \right )$, with $\lambda>0$.  For $\lambda
\ne 0$, let $E_\lambda\subset M_p$ be the direct sum of the 2-planes
with block $L_{\lambda}$.  Since $A$ is invariant under $\Hol$, so
is $E_\lambda$, and so is the orthogonal complement of $E_\lambda$.

By de Rham's holonomy theorem (see e.g\. \cite{B}), $M$ is then the
product of spaces $M_1$ and $M_2$ with tangent spaces $E_\lambda$ and
$E_\lambda^\perp$.  Since $M$ is irreducible, one of these spaces must
be trivial. Since $A \ne 0$, we can find a $\lambda$ such that the
vanishing subspace is $E_\lambda^\perp$, so $E_\lambda$ is all of
$M_p$. But then $A$ is $\lambda$ times a complex structure on $M_p$.
Furthermore, since $\Hol$ acts trivially on $A$, $\lambda^{-1}A$ can be
extended to a global integrable complex structure.  $M$ is thus a
complex symmetric space, and $\Hol$ is a subgroup of $U(k)$, with
$2k=n$.

Since $A \in \ker R^M$, and since $R^M$ is a self-adjoint map from
$\so(TM)$ to itself, $A$ is orthogonal to the image of $R^M$.  Thus,
for any $B\in\Lambda^2(TM)$,
$0=\<A,R^M(B)\>=\frac12\operatorname{tr}(A^t\circ R^M(B))$. But $A$ is
a (multiple of a) complex structure, so that each $R^M(B)\in\frak
u(k)$ has complex trace 0, and the holonomy of $M$ lies in $SU(k)$,
not just in $U(k)$. The determinant of the tangent bundle (viewed as a
complex $k$-dimensional vector bundle) is trivial, with a flat
connection.  In other words, $M$ is a Calabi-Yau manifold.

Calabi-Yau manifolds have Ricci-flat metrics (see Ch.9, Theorem 4.6 in
\cite{K-N}).  However, Ricci-flat homogeneous spaces (and in particular
Ricci-flat symmetric spaces) are flat (see section 7.61
in \cite{B}).  Since $M$ is simply connected, it must equal $\c^k$.
However, $\c^k = \r^k \times \r^k$ is reducible, and we have a
contradiction. \qed
\enddemo

\remark{Remark} The simple connectivity assumption in the previous
theorem was necessary to ensure that $\Hol$ is connected.  However,
a version of the theorem holds in the non-simply connected case. If $M$ is
not simply connected, we work on the universal cover $\tilde M$.  If
$\tilde M$ is irreducible, then Condition A holds on $\tilde M$, and 
therefore holds on $M$. 

\endremark

\remark { Remark} Condition A is considerably weaker than the
nonexistence of a parallel 2-form, which is a common condition applied
to symmetric spaces.  We have just seen that a parallel section of
$\so(TM)$ (or equivalently a parallel 2-form) on an irreducible
symmetric space must be a complex structure, but many irreducible
symmetric spaces are indeed complex, and their complex structures are
indeed parallel. Condition A allows such structures, just not in the
kernel of $R^M$.\endremark

To understand the general case, we will examine first how Condition A
behaves under products with irreducible symmetric spaces:

\proclaim{Lemma 4.2} Let $M\ne \r$ be a simply connected irreducible
symmetric space and $N$ a simply connected symmetric space, not
necessarily irreducible. Then $M\times N$ satisfies condition A if and
only if $N$ does.
\endproclaim

\demo{Proof\/} The curvature of $M\times N$ is the direct sum of the curvature
of $M$ and the curvature of $N$, so the kernel of $R^{M\times N}$ at a
point $(p,q)$ is the kernel of $R^{M}$ at $p$ plus the kernel of
$R^{N}$ at $q$ plus the linear subspace spanned by the wedge of $M_p$
and $N_q$.  Similarly, the image of $R^{M\times N}$ is the direct sum
of the images of $R^M$ and $R^N$ under the diagonal inclusion
$\so(M_p)\oplus\so(N_q)\subset\so(M_p\oplus N_q)$.

It is clear then that $[\ker R^M,\im R^{M\times N}]=[\ker R^M, \im
R^M]$ which spans $\ker R^M$. Similarly, $[\ker R^N,\im R^{M\times
N}]=[\ker R^N, \im R^N]$.  Since $M$ is irreducible and $M \ne \r$, it
is easy to see that $\left\{R^M(\so(M_p))M_p\right\}$ spans all of $M_p$.
Thus for any $x\in M_p$ we can find $a_i\in \so(M_p)$ and $x_i\in M_p$
with $\sum_iR^M(a_i)x_i=x$. This implies that
$$
x\wedge y=\sum_iR^M(a_i)x_i\w y=-\sum_i[x_i\wedge y, R^M(a_i)], \tag{4.1}
$$
and therefore $\span (M_p\wedge N_q)=\span([M_p\w N_q,\im R^M])$.  It
follows from this that the span of $[\ker R^{M\times N},\im R^{M\times
N}]$ equals  $\ker R^{M\times N}$ if and only if $\span([\ker R^N, \im
R^N])$ equals $\ker R^N$, giving us the conclusion of the lemma.
\qed
\enddemo

\proclaim{Theorem 1.1} Let $M$ be a simply connected symmetric space
with de Rham decomposition given by $M=\r^{n_0}\times M_1\times\cdots\times
M_l$. Then $M$ satisfies condition A if and only if $n_0\le
1$. Equivalently, $M$ satisfies condition A if and only if the
dimension of the center
of $G$ does not exceed 1.
\endproclaim

\demo{Proof\/} From the last two lemmas it is clear that $M_1\times\cdots\times
M_l$ satisfies condition A. On the other hand, the condition holds on
$\r^{n_0}$ precisely when $n_0\le1$. Thus, the first
part of the theorem follows trivially from Lemma 4.2. The statement
about the center follows, since the center of $G$ generates
the flat factors in the de Rham decomposition.
\qed
\enddemo

\head 5.\quad Examples and applications \endhead

\subhead A. Bundles over spheres \endsubhead

The sphere $S^n$ is a symmetric space with $G_0=SO(n+1)$, $H_0=SO(n)$,
$G=Spin(n+1)$ and $H=Spin(n)$. By the main theorem,
a parallel vector bundle over $S^n$ corresponds to a representation of
$Spin(n)$, and is thus the Whitney sum of vector bundles corresponding
to irreducible representations of $Spin(n)$.  Our analysis simplifies to
listing the irreducible representations of $Spin(n)$.

If $n=2$, then $SO(2)$ and $Spin(2)$ are both circle groups.  We
identify $SO(2)$ with $\r/(2\pi \z)$ and $Spin(2)$ with $\r/(4\pi \z)$,
so that the Lie algebras of $SO(2)$ and $Spin(2)$ are naturally identified.
Besides the trivial 1-dimensional representation, the irreducible
representations of $Spin(2)$ are all 2-dimensional, namely
$$ \rho_k(x) = \left ( \matrix \cos(kx/2) & \sin(kx/2)
\\ -\sin(kx/2) & \cos(kx/2) \endmatrix \right ).  \tag{5.1}
$$
The corresponding vector bundle is a plane bundle of Euler class $k$, or
equivalently a complex line bundle of first Chern class $k$.  Important
special cases are the Hopf fibrations $k=\pm 1$ and the tangent
bundle $k=2$.

Before analyzing the higher dimensional cases, we recall some general
concepts and facts about orthogonal representations, cf\. also
\cite{BD}. If $V$ is an inner product space, then an irreducible
orthogonal representation $\rho:H\to SO(V)$ is either unitary for some
complex structure, or else its complexification is irreducible, and
the latter is of {\sl real type\/}; i.e., it admits a {\sl real
structure\/}, or in other words, a complex conjugate linear map
$J:\Bbb C\otimes V\to \Bbb C\otimes V$ that squares to the identity
and commutes with the action of $H$. Its restriction to the +1
eigenspace of $J$ is the original real representation. A unitary
representation is said to be of {\sl quaternionic type\/} if it admits
a {\sl quaternionic structure\/}, or in other words, a complex
conjugate linear map that squares to minus the identity and commutes
with $H$.

If $n=3$, then $Spin(3)=SU(2)$, and the irreducible representations
are well known, see for example \cite{BD}.  The fundamental
representation, corresponding to matrix multiplication, has complex
dimension 2, and has a quaternionic structure given by $J(z_1,z_2) =
(-\bar z_2,\bar z_1)$ for $(z_1,z_2)\in\Bbb C^2$.  All other complex
irreducible representations are obtained from symmetric products of
the fundamental representation.  Explicitly, let $V_k$ denote the
space of homogeneous polynomials $p$ of degree $k$ in two complex
variables,
$$
p(z_1,z_2)=\sum_{r=0}^k a_rz_1^rz_2^{k-r},\tag{5.2}
$$
of complex dimension $k+1$. The action of $SU(2)$ on $V_k$ is given by
$$
gp(z)=p(gz),\qquad g\in SU(2),\quad z\in\Bbb C^2,\quad p\in V_k.\tag{5.3}
$$
In particular, the structure map $J$ on $\Bbb C^2=V_1$ induces a
structure map $J$ on $V_k$ satisfying $J^2p=(-1)^kp$. Thus, the
odd-dimensional representations have an underlying real structure,
while the (complex) even-dimensional representations have a
quaternionic structure.  Complex irreducible representations of odd
dimension $2 \ell+1$ therefore correspond to real irreducible
representations of dimension $2 \ell +1$, while complex irreducible
representations of dimension $2\ell$ correspond to real irreducible
representations of dimension $4\ell$.  The odd-dimensional irreducible
representations also are representations of $SO(3)$ (since $-1\in
SU(2)$ acts trivially on $V_{2l}$, so that the action of $SU(2)$
induces one of $SO(3)=SU(2)/\pm1$), whereas the even-dimensional
irreducible representations are not.

Over $S^3$, therefore, there exists a unique irreducible parallel
unitary bundle in each (complex) dimension.  Irreducible Riemannian
parallel vector bundles exist in odd dimensions, and in dimensions
divisible by 4.  The first few examples are:
\item {1.} The trivial bundle, of real dimension 1.
\item {2.} The tangent bundle, of dimension 3.  This corresponds to the
fundamental representation of $SO(3)$.  $\Lambda^2(TM)$ is also 3-dimensional,
and is isomorphic to $TM$ via the Hodge star operator.
\item {3.} The spinor bundle, of complex dimension 2 and real dimension 4,
corresponding to the fundamental representation of $SU(2)$.
\item {4.} The quadrupole bundle of traceless symmetric endomorphisms
from $TM$ to itself, of real dimension 5.

If $n=4$, then $Spin(n)=SU(2)\times SU(2)$.  The irreducible representations
are tensor products of the representations of the first and second factors.
If we take the $k_1+1$-dimensional representation of the first factor and
the $k_2+1$-dimensional representations of the second, we obtain a
representation of complex dimension $(k_1+1)(k_2+1)$ with a structure map
$J$ that squares to $(-1)^{k_1+k_2}$.  We label this representation
$(k_1,k_2)$.  The complex dimension of $(k_1,k_2)$ is
always $(k_1+1)(k_2+1)$, while the real dimension is $(k_1+1)(k_2+1)$ if
$k_1+k_2$ is even and  $2(k_1+1)(k_2+1)$ if $k_1+k_2$ is odd.  Once again,
we see that real irreducible representations come only in odd dimensions
and in dimensions that are multiples of 4.

The first few examples of irreducible parallel Riemannian bundles are:
\item {1.} The trivial bundle, of dimension 1, corresponding to the
representation $(0,0)$.
\item {2.} The bundle $\Lambda^+$ of self-dual 2-forms on $S^4$, of
real dimension 3, corresponding to the representation $(2,0)$.
\item {3.} The bundle $\Lambda^-$ of anti-self-dual 2-forms on $S^4$, of
real dimension 3, corresponding to the representation $(0,2)$.
\item {4.} The tangent bundle, of real dimension 4, corresponding to
the representation $(1,1)$.
\item {5.} The bundle $S^+$ of positive spinors, of complex dimension 2
and real dimension 4, originating from the representation
$(1,0)$.  The parallel connection is the standard round self-dual
$SU(2)$ instanton of second Chern class $-1$.
If we identify $S^4 = {\Bbb H}{\Bbb P}^1$, this is
the tautological quaternionic line bundle, and the unit sphere bundle is
a Hopf fibration $S^7 \to S^4$.
\item {6.} The bundle $S^-$ of negative spinors,
corresponding to the representation $(0,1)$.  The
parallel connection, with structure group $SU(2)$, is the round anti-self-dual
anti-instanton, of second Chern class $1$.
\item {7.} There are bundles of real dimension 5, corresponding to the
representations $(4,0)$ and $(0,4)$.  These are the traceless symmetric
endomorphisms of $\Lambda^+$ and $\Lambda^-$, respectively.

Combining these low-dimensional irreducible representations, we obtain
a list of all the rank-4 Riemannian vector bundles over $S^4$ with
parallel curvature.

\proclaim{Proposition 5.1} There are exactly six rank-4 Riemannian vector
bundles over $S^4$ with parallel curvature.  Three of these, namely
$S^+$, $S^-$ and the trivial bundle, admit a complex structure.  The
corresponding rank-2 complex vector bundles have second Chern class $-1$,
$+1$ and $0$, respectively.  The remaining three Riemannian vector bundles,
namely $TM$, $\r \oplus \Lambda^+$ and $\r \oplus \Lambda^-$, are not complex.
\endproclaim

\demo{Proof} These examples are constructed from representations of $Spin(4)$,
and there are no other real 4-dimensional representations.  All that
remains is to show that these bundles are distinct.  This can be seen from
their characteristic classes. $TM$ and $S^\pm$ all have nontrivial Euler
classes (2 and $\mp1$ respectively), and their Pontryagin classes are 0,
$\pm2$. The bundles $\r^4$, $\r\oplus\Lambda^+$ and
$\r\oplus\Lambda^-$ all have trivial Euler class, but their Pontryagin
classes are all distinct (0 and $\pm4$, respectively).
\qed
\enddemo

If $n>4$, then $SO(n)$ is simple, and $spin(n)$ is a simple Lie Algebra.
The smallest nontrivial representation is the fundamental representation
of $SO(n)$.  Hence there are no nontrivial parallel bundles of dimension
less than $n$, and the only nontrivial parallel bundle of dimension $n$
is the tangent bundle.  Other noteworthy representations, and their
corresponding bundles, are
\item {1.} The $k$-th exterior power of the fundamental representation
of $SO(n)$ has dimension $n \choose k$, and corresponds to the bundle
of $k$-forms on $S^n$.
\item {2.} The symmetric product of the fundamental representation of $SO(n)$
with itself splits as a trivial factor plus an $[n(n+1)-2]/2$ irreducible
representation.  The latter corresponds to the traceless symmetric
endomorphisms of $TM$.
\item {3.} The fundamental representation of $Spin(n)$ has complex
dimension $2^{n/2}$ if $n$ is even and $2^{(n+1)/2}$ if $n$ is odd.
If $n$ is even the spin bundle splits as the sum of positive and
negative spinor bundles, each of complex dimension $2^{(n-2)/2}$.

\subhead B. Bundles over complex projective spaces \endsubhead

The complex projective space $M=\CP^n$ is naturally the quotient of $U(n+1)$
by $U(1) \times U(n)$.  However, this quotient does not describe the
structure of $\CP^n$ as a symmetric space, since the center of $U(n+1)$
is contained in $U(1)\times U(n)$.  We must take the quotient of both
$U(n+1)$ and $U(1)\times U(n)$ by $U(1)$.  The question is, which $U(1)$?

The most natural choice of $U(1)$ is the center $Z$ of $U(n+1)$.  Let
$G_0 = U(n+1)/Z$, and let $H_0=(U(1)\times U(n))/Z \sim U(n)$.  $G_0/H_0$
is indeed $\CP^n$, $H_0$ is the holonomy group of $\CP^n$, and $G_0$ is
the bundle of complex orthonormal frames of $TM$.

$G_0$ is not simply connected, however.  Its universal cover is
$G=SU(n+1)$, and the covering map is $n\!+\!1$-to-1.  $G$ can be
viewed as the quotient of $U(n+1)$ by $Z'=U(1) \times I_n$, where
$I_n$ is the $n\times n$ identity matrix, while $H= (U(1)\times
U(n))/Z'\sim U(n)$.  Although $H$ and $H_0$ are both isomorphic to
$U(n)$ as groups, $H$ is the $n+1$-fold cover of $H_0$, with
covering map $A \to (\det A)A$.

$G$ itself can be viewed as a bundle of frames, just not of $TM$.  $G$ is
the bundle of complex orthonormal frames of $L \otimes_\c TM$, where $L$
is the tautological complex line bundle over $\CP^n$.  Equivalently,
$L \otimes_\c TM$ is the orthogonal complement to $L$ in the trivial
$\c^{n+1}$ bundle over $\CP^n$.

There are infinitely many nontrivial complex 1-dimensional (or real
2-dimensional) representations of $H=U(n)$, namely all powers of the
determinant.  The determinant representation itself yields a line
bundle of Chern class $+1$.  This is the determinant bundle of $L
\otimes_\c TM$, and is also isomorphic to $L^*$, the complex dual of
$L$. The $k$-th power of the determinant yields $(L^*)^k$ if $k>0$, or
$L^{-k}$ if $k<0$.  (For $n=1$, this replicates our previous analysis
of $M=S^2$.)

Note that the determinant representation of $H_0$ yields the
determinant of $TM$ (viewed as a complex vector bundle), which has
Chern class $n+1$, and all complex line bundles constructed from
representations of $H_0$ have first Chern class divisible by
$n+1$. This demonstrates once again that to get all parallel bundles
we must consider all representations of $H$, not just representations
of $H_0$.

After the trivial bundle, the determinant, and powers of the
determinant, the smallest unitary representations of $U(n)$ are the
fundamental representation and its dual, both of complex dimension
$n$.  There are an infinite number of complex representations of this
dimension, obtained by multiplying the trivial representation (or its
dual) by arbitrary powers of the determinant.  The smallest
nontrivial real (i.e. not complex) irreducible representation is the
adjoint representation, of real dimension $n^2-1$.

\newpage
\subhead C. Lie Groups \endsubhead

Every connected Lie group $K$ which is the product of a compact group
and a vector group is a symmetric space.  If $K$ is simply connected,
we can take $G=K\times K$, and let $H\sim K$ be the diagonal subgroup
of $G$. Parallel bundles over $K$ thus correspond to representations
of $K$ itself.  The tangent bundle of $K$ corresponds to the adjoint
representation, as the tangent space of $K$ at the origin is naturally
isomorphic to the Lie algebra of $K$. It follows from Corollary 5.3
below that if $K$ has no Euclidean factor, then the unit tangent
bundle of $K$ admits a metric of constant positive scalar curvature.

Applying this construction to $K=SU(2)$, we recover our previous results
about $M=S^3$.

\subhead D. Sphere bundles and constant scalar curvature \endsubhead

Let $E$ denote a rank $k$ vector bundle with parallel connection over
$M=G/H$. It follows from the main theorem that if $G$ admits a
bi-invariant metric, then $E$ and the unit sphere bundle $E^1$ of $E$
admit complete metrics of nonnegative sectional curvature: Endow $\Bbb
R^k$ with any $O(k)$ invariant metric. Then the diagonal action of $H$
on the metric product $G\times\Bbb R^k$ is by isometries, and there is
a unique metric on the quotient $E=G\times_\rho\Bbb R^k$ (respectively
$G\times_\rho S^{k-1}$) for which the projection $G\times \Bbb R^k\to
E$ (resp\. $G\times S^{k-1}\to E^1$) becomes a Riemannian submersion,
see \cite{C}. Moreover, if $H$ acts transitively on the unit sphere in
$\Bbb R^k$, then $E^1$ is a homogeneous space under the isometric
action $g[(g',u)]:=[(gg',u)]$, $g\in G$, and in particular, $E^1$ has
constant scalar curvature. It is noteworthy that there are vector
bundles over the sphere that admit nonnegatively curved metrics but no
parallel connections: One of these, in fact, is a rank 4 bundle over
$S^4$ whose unit sphere bundle is a Milnor exotic 7-sphere, see
\cite{G-M} and
\cite{R}.

If $H$ does not act transitively on $S^{k-1}$, then the resulting
metric on the unit sphere bundle is not, in general,
homogeneous. Nevertheless, when the action of $H$ is irreducible, then
there exists a metric of {\sl constant scalar curvature\/} on $E^1$:
Decompose $TE=\Cal H\oplus\Cal V$ as a sum of subbundles $\Cal H$ (the
horizontal distribution defining the connection) and $\Cal V$ (the
vertical distribution tangent to the fibers of $\pi:E\to M$) over
$E$. The inner product on $\Cal H$ is defined by pulling back the one
on $TM$ via $\pi_*:\Cal H\to TM$. The one on $\Cal V$ is induced by
the $O(k)$ invariant metric on $\Bbb R^k$. Finally, declare both
distributions to be mutually orthogonal. This defines a {\sl
connection metric\/} on $E$, and the projection $\pi:E\to M$ is a
Riemannian submersion with totally geodesic fibers. It can be shown
\cite{S-W} that if the sectional curvature $K_M$ of $M$ is strictly
positive (i.e., if $M$ is a compact rank one symmetric space), then
$K_E\ge0$, provided the fiber metric is sufficiently curved.
\proclaim{Proposition 5.2} The connection metric on $E$ induces a constant
scalar curvature metric on the unit sphere bundle $E^1$.
\endproclaim

\demo{Proof} Let $s_E$, $s_M$, and $s_F$ denote the scalar curvatures of $E$,
$M$, and the fiber $F$ respectively. Since $\pi:E\to M$ has totally
geodesic fibers,
$$
s_E=s_m\circ\pi+s_F-|A|^2,\qquad |A|^2:=\sum_{i,j}|A_{x_i}x_j|^2,\tag5.4
$$
cf\. \cite{B, 9.37}. Here $A$ denotes the O'Neill tensor of $\pi$,
$A_XY=\frac12[X,Y]^v$ for horizontal $X$ and $Y$, and $x_i$ is (the
horizontal lift of) an orthonormal basis of $TM$. On the other hand, for unit
$u\in E$ and $r>0$,
$$
|A|^2(ru)=\frac14G^2(r)\sum_{i,j}|R^E(x_i,x_j)u|^2,\tag5.5
$$
where $dr^2+G^2(r)d\sigma^2$ is the fiber metric in polar coordinates, see
\cite{S-W}. It therefore suffices to check that the map $C:E^1\to\Bbb R^+$
given by $C(u)=\sum|R(x_i,x_j)u|^2$ is constant. Observe that $C(u)=\<\tilde
Cu,u\>$, where $\tilde C:E\rightarrow E$ is defined as
$$
\tilde C=-\sum_{i,j}R^E(x_i,x_j)\circ R^E(x_i, x_j).\tag5.6
$$
But $\tilde C$ is parallel since $R^E$ is, and therefore commutes with
the action of the holonomy group of $E$. The latter is by assumption
irreducible, so that by Schur's lemma, $\tilde C$ is a multiple of the
identity. This implies that $C$ is constant as claimed.
\enddemo

\proclaim{Corollary 5.3} Let $M$ be an irreducible symmetric space with
nonnegative
sectional curvature. Then the unit tangent bundle of $M$ admits a metric of
nonnegative sectional curvature and constant scalar curvature.
\endproclaim

\demo{Proof} The first part of the statement is essentially due to Cheeger
\cite{C}. The claim about the scalar curvature is an immediate consequence of
Proposition 5.2 together with the irreducibility assumption.
\enddemo

\noindent{\bf Acknowledgements.} We thank Sharad Agnihotri and
Jan Segert for helpful discussions.  The work of L.S. is partially
supported by NSF grant DMS-9626698.

\enddocument
\bye